\documentclass[12pt]{amsart}
\usepackage[latin1]{inputenc}
\usepackage{mathptmx}
\usepackage{amscd}
\usepackage{amssymb}
\newtheorem{theorem}{Theorem}
\newtheorem{lemma}[theorem]{Lemma}
\newtheorem{corollary}[theorem]{Corollary}
\newtheorem{proposition}[theorem]{Proposition}
\newtheorem{conjecture}[theorem]{Conjecture}

\theoremstyle{definition}
\newtheorem{remark}[theorem]{Remark}
\newtheorem{definition}[theorem]{Definition}

\newtheorem{examples}[theorem]{Examples}

\DeclareMathOperator{\Id}{Id}

\let\epsilon=\varepsilon

\begin{document}
	
	\title{Amount algebras}
	
	\author{Peyman Nasehpour}
	
	\address{Department of Engineering Sciences, Golpayegan University of Technology, Golpayegan, Isfahan Province, Iran}
	
	\email{nasehpour@gut.ac.ir, nasehpour@gmail.com}

	\begin{abstract}
		In this paper, as a generalization to content algebras, we introduce amount algebras. Similar to the Anderson-Badawi $\omega_{R[X]}(I[X])=\omega_R(I)$ conjecture, we prove that under some conditions, the formula $\omega_B(I^{\epsilon})=\omega_R(I)$ holds for some amount $R$-algebras $B$ and some ideals $I$ of $R$, where $\omega_R(I)$ is the smallest positive integer $n$ that the ideal $I$ of $R$ is $n$-absorbing. A corollary to the mentioned formula is that if, for example, $R$ is a Pr\"{u}fer domain or a torsion-free valuation ring and $I$ is a radical ideal of $R$, then $\omega_{R[][X]]}(I[[X]])=\omega_R(I)$. 
	\end{abstract}
	
	\maketitle
	
	\section{Introduction}
	
	In this paper, all rings are commutative with identity and all algebras are unitary \cite{Matsumura1989}. Let us recall that a proper ideal $I$ of a ring $R$ is an $n$-absorbing ideal of $R$, if whenever $x_1 \cdots x_{n+1} \in I$ for $x_1, \ldots, x_{n+1} \in R$, then there are $n$ of the $x_i$'s whose product is in $I$. Anderson and Badawi \cite{AndersonBadawi2011} conjectured that \[\omega_{R[X]}(I[X])=\omega_R(I) \qquad \text{(Anderson-Badawi $\omega$ Conjecture)}\] for each ideal $I$ of an arbitrary ring $R$, where \[\omega_R(I)= \min \{n \colon \text{ $I$ is an $n$-absorbing ideal of $R$}\}.\] In this direction, the author proved that if $R$ is a Pr\"{u}fer domain, then for any content $R$-algebra $B$, $\omega_B(IB) = \omega_R(I)$ and since any polynomial ring $R[X]$ is a content $R$-algebra (see Hilfsatz von Dedekind-Mertens on p. 128 in \cite{Krull1968}), it is clear that the Anderson-Badawi $\omega$ conjecture is true if $R$ is a Pr\"{u}fer domain \cite[Corollary 11]{Nasehpour2016AB}. The main purpose of this paper is to prove that under some conditions the formula $\omega_{R[[X]]}(I[[X]])=\omega_R(I)$ holds as well. In fact, inspired by the recent papers of Epstein and Shapiro \cite{EpsteinShapiro2016} and Kang et al. \cite{KangParkToan2018}, we introduce {\em amount algebras} and show that under some conditions - that we are going to report in the upcoming passages - some formulas similar to  $\omega_{R[X]}(I[X])=\omega_R(I)$ holds in amount algebras and a corollary to these results is that under some conditions $\omega_{R[[X]]}(I[[X]])=\omega_R(I)$ is also true. Here is a brief sketch of the contents of our paper:
	
	In Definition \ref{amonutfunctionsdef}, we introduce the concept of amount functions as follows: 
	
	Let $R$ be a ring and $B$ an $R$-algebra. We say a function $A$ from $B$ to the set of ideals $\Id(R)$ of $R$ defined by $f \mapsto A_f$ is an amount function if the following properties hold for all $r\in R$ and $f,g\in B$:
	
	\begin{enumerate}
		\item $A$ preserves 0 and 1, i.e. $A_0 = (0)$ and $A_1 = R$.
		\item If $A_f =(0)$ then $f=0$. 
		\item $A$ is homogeneous, i.e. $A_{rf} = rA_f$.
		\item $A$ is submultiplicative, i.e. $A_{fg} \subseteq A_f A_g$.
	\end{enumerate}
	
	A general example for amount functions is the content function $c$ over a faithfully flat $R$-algebra $B$ with this additional property that $B$ as an $R$-module is content (check Theorem \ref{contentisamount1}). Other examples (see Examples \ref{examplesamountfunctions}) include the function $A$ defined on power series rings $R[[X]]$ by $A_f = (r_0,r_1,\ldots, r_n,\ldots),$ where $f= r_0+r_1 X + \cdots +  r_n X^n + \cdots$ is an element of $R[[X]]$ \cite{GilmerGramsParker1975}. On the other hand, for all $f,g \in R[[X]]$, we have the following {\em amount formulas}:
	\begin{itemize}
		\item $A_f ^{n+1} A_g = A_f ^n A_{fg}$, for some $n$, if $R$ is Noetherian and $n\in \mathbb N_0$ depending on $g$ is large enough \cite[Theorem 2.6]{EpsteinShapiro2016}.
		
		\item $A^2_f A_g = A_f A_{fg}$ or $A^2_g A_f = A_g A_{gf}$ if $D$ is a valuation ring \cite[Theorem 2.8]{KangParkToan2018}.
		
		\item $(A_f A_g)^2 = A_f A_g A_{fg}$ if $D$ is a Pr\"{u}fer domain \cite[Corollary 2.9]{KangParkToan2018}).
	\end{itemize}

	Inspired by the amount formulas mentioned in above, we define amount algebras (check Definition \ref{amountalgebras}) as follows:
	
	Let $R$ be a ring and $B$ an $R$-algebra. We say $B$ is an amount $R$-algebra if the following conditions hold:
	
	\begin{itemize}
		\item  There is an amount function $A$ from $B$ to $\Id(R)$ defined by $f \mapsto A_f$ with this property that for all $f,g \in B$, there are non-negative integers $m,n$ such that \[A_f ^m A_g ^n A_{fg} = A_f ^{m+1} A_g ^{n+1}.\]
		
		\item There is a function $\epsilon$ from $\Id(R)$ to $\Id(B)$ defined by $I \mapsto I^{\epsilon}$ with the following properties:
		
		\begin {enumerate}
		\item $A_f \subseteq I$ if and only if $f\in I^{\epsilon}$, for all $f\in B$ and $I\in \Id(R)$.
		\item $I^{\epsilon} \cap R = I$, for all $I\in \Id(R)$.
		\end {enumerate}
	\end{itemize}
	
	Let us recall that an ideal $I$ of a commutative ring $R$ is strongly $n$-absorbing if whenever \[I_1 \cdots I_{n+1} \subseteq  I\] for some ideals $I_1, \ldots, I_{n+1}$ of $R$, then there are $n$ of the $I_i$'s whose product is a subset of $I$.
	
	In Theorem \ref{ABconjecture1}, we prove that if $R$ is a ring such that any $n$-absorbing ideal $I$ of $R$ is strongly $n$-absorbing for any positive integer $n$, also $B$ is an amount $R$-algebra, and $B$ is Gaussian, then $\omega_B(I^{\epsilon}) = \omega_R(I)$. A corollary (see Corollary \ref{AB1cor}) to this is that if $I$ is an ideal of a Dedekind domain $D$, then \[\omega_{D[[X]]} (I[[X]]) = \omega_{D[X]} (I[X]) =  \omega_D (I).\] Note that an amount $R$-algebra $B$ is Gaussian if $A_{fg} = A_f A_g$ for all $f,g \in B$ (check Definition \ref{Gaussianamountalgebra}).
	
	Also in Theorem \ref{ABconjecture2}, we show that if $R$ is a ring such that any $n$-absorbing ideal $I$ of $R$ is strongly $n$-absorbing for any positive integer $n$, $B$ is an amount $R$-algebra, and $I$ is a radical ideal of $R$, then $\omega_B(I^{\epsilon}) = \omega_R(I)$. A corollary (see Corollary \ref{AB2cor1} and Corollary \ref{AB2cor2}) to this result is that if $I$ is a radical ideal of a ring $R$, and either $R$ is a torsion-free Noetherian  ring, or $D$ is a Pr\"{u}fer domain, or a torsion-free valuation ring, then \[\omega_{D[[X]]} (I[[X]]) = \omega_{D[X]} (I[X]) =  \omega_D (I).\] We end our paper by conjecturing that if $I$ is an ideal of a ring $R$, then \[\omega_{R[[X]]} (I[[X]]) = \omega_R (I).\]

	\section{Amount Algebras}\label{sec:amount}
	
	We begin this section by introducing the amount functions.
	
	\begin{definition}[Amount functions] \label{amonutfunctionsdef} Let $R$ be a ring and $B$ an $R$-algebra. We say a function $A$ from $B$ to the set of ideals $\Id(R)$ of $R$ is an amount function if the following properties hold for all $r\in R$ and $f,g\in B$:
		
		\begin{enumerate}
			\item $A$ preserves 0 and 1, i.e. $A_0 = (0)$ and $A_1 = R$.
			\item If $A_f =(0)$ then $f=0$. 
			\item $A$ is homogeneous, i.e. $A_{sf} = sA_f$.
			\item $A$ is submultiplicative, i.e. $A_{fg} \subseteq A_f A_g$.
		\end{enumerate}
		
	\end{definition}

	\begin{examples} \label{examplesamountfunctions} In the following, we bring two important examples for amount functions:
		
		\begin{enumerate}
			\item Let $(\Gamma,+,0,<)$ be a totally ordered commutative additive monoid and $R$ be a ring. Let $f=r_1 X^{\alpha_1} + r_2 X^{\alpha_2} + \cdots + r_n X^{\alpha_n}$ be an element of the monoid ring $R[\Gamma]$. Define the content of $f$, denoted by $c(f)$, to be an ideal of $R$ generated by the coefficients of $f$, i.e. \[c(f) := (r_1,r_2,\ldots,r_n).\] It is easy to verify that $c: R[\Gamma] \longrightarrow \Id(R)$ is an amount function. Note that by $\Id(R)$, we mean the set of all ideals of the ring $R$.
			
			\item Let us recall that an element $x$ of a totally ordered semigroup $(\Gamma,+, <)$ is finitely decomposable if there are only finitely many pairs $(y_i,z_i)$ of elements of $\Gamma$ such that $x=y_i + z_i$. Now, let $(\Gamma, +, 0, <)$ be a totally ordered additive commutative monoid. Assume that 0 is the least element of $\Gamma$ and that each element of $\Gamma$ is finitely decomposable (for example, let $\Gamma = \bigoplus \mathbb N_0$). Let $R$ be a ring and $R[[\Gamma]]$ be the set of all functions $f: \Gamma \rightarrow R$. Let $f$ and $g$ be arbitrary elements of $R[[\Gamma]]$ and define their addition and multiplication as follows: \[(f+g)(x) = f(x)+g(x), ~ (fg)(x) = \sum_{y+z=x} f(y)g(z).\] It is straightforward to see that $R[[\Gamma]]$ is an $R$-algebra \cite{GilmerGramsParker1975}. For each $f\in R[[\Gamma]]$, define $A_f$ to be an ideal of $R$ generated by all $f(s)$, i.e. coefficients of $f$. It is easy to see that the function $A$ from $R[[\Gamma]]$ to $\Id(R)$ defined by $A \mapsto A_f$ is an amount function. For instance, for an element $f=s_0 + s_1 X + \cdots + s_n X^n + \cdots $ in $R[[X]]$, \[A_f = (s_0, s_1, \ldots, s_n, \ldots).\]
		\end{enumerate}

	\end{examples}
	
	Let us recall that if $B$ is an $R$-algebra. The content function $c: B \rightarrow \Id(R)$ is defined by \[c(f) = \bigcap \{I\in \Id(R): f\in IB\},\] where by $IB$, we mean the extension of the $R$-ideal $I$ in $B$. By definition, $B$ as an $R$-module is content if $f\in c(f)B$ for all $f\in B$ \cite{OhmRush1973}.
	
	\begin{theorem} \label{contentisamount1}
		Let $B$ be an $R$-algebra and a content $R$-module. The content function $c$ is an amount function if and only if $B$ is a faithfully flat $R$-module.
		
		\begin{proof} Let $B$ be an $R$-algebra. It is clear that $c(0)= (0)$. If $B$ is a content $R$-module, then $f\in c(f)B$ and $g\in c(g)B$, for arbitrary elements $f$ and $g$ in $B$ and so, $fg\in c(f)c(g)B$. This implies that $c(fg) \subseteq c(f)c(g)$ (see Proposition 1.1 in \cite{Rush1978}). On the other hand, $B$ is flat if and only if $c(rf) = rc(f)$ for all $r\in R$ and $f\in B$ \cite[Corollary 1.6]{OhmRush1973}. Also, according to Corollary 1.6 and the Statement 6.1(a) in \cite{OhmRush1973} and Proposition 1.1 in \cite{Rush1978}, if $B$ is a content and flat $R$-module, then $B$ is faithfully flat if and only if $c(1) = R$ and the proof is complete. \end{proof}
	\end{theorem}
	
	\begin{remark}
		If an $R$-algebra $B$ as a module is content, then $c(f)$ is finitely generated for all $f\in B$ \cite[\S1]{OhmRush1973}. Now, let $R[[\Gamma]]$ be as the $R$-algebra defined in Examples \ref{examplesamountfunctions}. It is clear that for $f\in R[[\Gamma]]$, the ideal $A_f$ is not necessarily finitely generated. 
	\end{remark}
	
	The proof of the following is straightforward:
	
	\begin{proposition}
		Let $B$ be an $R$-algebra and $A$ an amount function from $B$ to $\Id(R)$. Then the following statements hold:
		
		\begin{enumerate}
			\item $A_r = (r)$ for all $r\in R$. In particular in Definition \ref{amonutfunctionsdef}, the condition $A_0 = (0)$ is superfluous.
			\item The equality $A_f A_g = (0)$ implies $fg=0$ for all $f,g \in B$.
		\end{enumerate}
	\end{proposition}
	
	Now we define amount algebras:
	
	\begin{definition} \label{amountalgebras}
		Let $R$ be a ring and $B$ an $R$-algebra. We say $B$ is an amount $R$-algebra if the following conditions hold:
		
		\begin{enumerate}
			\item  There is an amount function $A$ from $B$ to $\Id(R)$ defined by $f \mapsto A_f$ with this property that for all $f,g \in B$, there are non-negative integers $m,n$ such that \[A_f ^m A_g ^n A_{fg} = A_f ^{m+1} A_g ^{n+1}\qquad \text{(The Amount Formula).} \label{Amountformula} \]
			
			\item There is a function $\epsilon$ from $\Id(R)$ to $\Id(B)$ defined by $I \mapsto I^{\epsilon}$ with the following properties:
			
			\begin{enumerate}
			\item $A_f \subseteq I$ if and only if $f\in I^{\epsilon}$, for all $f\in B$ and $I\in \Id(R)$.
			\item $I^{\epsilon} \cap R = I$, for all $I\in \Id(R)$.
			\end{enumerate}
		\end{enumerate}
	\end{definition}

	\begin{proposition}
		Let $B$ be an amount $R$-algebra. Then the following statements hold:
		
		\begin{enumerate}
			\item $f\in A_f^{\epsilon}$ for all $f\in B$.
			\item $I \subseteq J$ if and only if $I^{\epsilon} \subseteq J^{\epsilon}$ for all ideals $I$ and $J$ of $R$.
		\end{enumerate}
		
		\begin{proof} (1): Since $A_f \subseteq A_f$, by definition, $f\in A_f^{\epsilon}$.
			
			(2): Assume that $I \subseteq J$ and let $f\in I^{\epsilon}$. By definition, $A_f \subseteq I$. So, $A_f \subseteq J$. This implies that $f\in J^{\epsilon}$. On the other hand, if $I^{\epsilon} \subseteq J^{\epsilon}$, then $I^{\epsilon} \cap R \subseteq J^{\epsilon} \cap R$ which is equivalent to say that $I \subseteq J$.\end{proof}
		
	\end{proposition}
	
	Let ue recall that if $I$ and $J$ are ideals of a ring $R$ then $J$ is a reduction of $I$ if $J\subseteq I$ and $JI^k = I^{k+1}$ for some positive integer $k$ \cite[Definition 1]{NorthcottRees1954}.
	
	\begin{lemma}
		
		\label{amountreduction}
		Let $B$ be an amount $R$-algebra. Then $A_{fg}$ is a reduction of $A_f A_g$ for all $f,g\in B$.
		
		\begin{proof}  Let $f,g\in B$. Then by definition, there are non-negative integers $m,n$ such that $A_f ^m A_g ^n A_{fg} = A_f ^{m+1} A_g ^{n+1}.$ Let $k=1+\max\{m,n\}$. So, $A_{fg} (A_f A_g) ^{k} = (A_f A_g) ^{k+1}$. Clearly, $k$ is a positive integer and $A_{fg} \subseteq A_f A_g$. Hence, $A_{fg}$ is a reduction of $A_f A_g$ and the proof is complete.\end{proof}
	\end{lemma}

\begin{theorem}
	Let $B$ be an amount $R$-algebra. Then $A_f A_g \subseteq \sqrt{A_{fg}}$ for all $f,g\in B$.
	
	\begin{proof} Let $P$ be a prime ideal of $R$ containing $A_{fg}$. By Lemma \ref{amountreduction}, $A_{fg}$ is a reduction of $A_f A_g$. So, $A_{fg} (A_f A_g) ^{k} = (A_f A_g) ^{k+1}$ for some positive integer $k$. This implies that $P$ contains $A_f A_g$. Hence, $\displaystyle A_f A_g \subseteq \bigcap_{P \supseteq A_{fg}} P = \sqrt{A_{fg}}$. This completes the proof. \end{proof}
\end{theorem}
	
	Let $B$ be an $R$-algebra such that as an $R$-module, it is content and faithfully flat. Then, $B$ is called to be a content $R$-algebra \cite[\S6]{OhmRush1973} if for all $f,g\in B$, there is a non-negative integer $n$ such that the Dedekind-Mertens formula $c(f)^{n+1} c(g) = c(f)^n c(fg)$ holds. 
	
	\begin{theorem}
		\label{contentisamount2}
		Let $B$ be a content $R$-algebra. Then $B$ is an amount $R$-algebra.
		
		\begin{proof}  Assume that $B$ is a content $R$-algebra. By Theorem \ref{contentisamount1}, $c(f)$ is an amount function. Obviously, the Dedekind-Mertens formula is a kind of the amount formula given in Definition \ref{amountalgebras}. Now, define $I^{\epsilon} = IB$. Clearly, $c(f) \subseteq I$ if and only if $f\in IB$ for all $f\in B$ and $I\in \Id(R)$, since $c(f)$ is the smallest ideal satisfying the condition $f\in IB$ \cite[\S1]{OhmRush1973}. Finally, it is clear that $I \subseteq IB \cap R$. Now, let $r\in IB \cap R$. So, $c(r) \subseteq I$. But $c(r) = (r)$ for all $r\in R$. Therefore, $r\in I$. Hence, $ IB \cap R \subseteq I$. From all we said, we conclude that $B$ is an amount $R$-algebra and the proof is complete. \end{proof}
	\end{theorem}
	
	Let $(\Gamma,+,0,<)$ be a totally ordered commutative additive monoid and $R$ be a ring. Northcott \cite{Northcott1959} has proved that $R[\Gamma]$ is a content $R$-algebra. Consequently, we have the following corollary:

	\begin{corollary}
		If $(\Gamma,+,0,<)$ is a totally ordered commutative additive monoid and $R$ is a ring, then the monoid ring $R[\Gamma]$ is an amount $R$-algebra.
	\end{corollary}
	
	\begin{remark}[More examples for amount algebras] Let $R$ be a ring and $X$ an indeterminate over $R$. Define $A_f$ to be the $R$-ideal generated by the coefficients of $f$ in the power series ring $R[[X]]$ and set $I^{\epsilon} = I[[X]]$. Note that $I[[X]]$ is not in general equal to $I\cdot R[[X]]$ \cite[Proposition 1]{GilmerHeinzer1968}). Now, it is easy to verify that all the properties necessary for $R[[X]]$ to be an amount $R$-algebra hold except the possibility of the amount formula given in Definition \ref{amountalgebras}. However, $R[[X]]$ is an amount $R$-algebra if $R$ is either Noetherian \cite[Theorem 2.6]{EpsteinShapiro2016}, or a Pr\"ufer domain \cite[Corollary 2.9]{KangParkToan2018}, or a valuation ring \cite[Theorem 2.8]{KangParkToan2018}.
		
	\end{remark}
	
	\begin{definition} \label{Armendarizamountalgebra}
		We say an amount $R$-algebra $B$ is Armendariz if $fg=0$ implies $A_f A_g =(0)$ for all $f,g \in B$, where $A$ is the amount function defined in Definition \ref{amonutfunctionsdef}. 
	\end{definition}
	
	Let us recall that a ring $R$ is reduced if $r^n = 0$ for some $n\in \mathbb N$ implies $r=0$ \cite[p. 3]{Matsumura1989}.
	
	\begin{theorem} \label{AmountisArmendariz}
		Let $R$ be a reduced ring and $B$ an amount $R$-algebra. Then $B$ is Armendariz. In particular, for all $f\in B$, we have the following: \[f\in Z_B(B) \implies fr=0 \hbox{~for some $r$ in $R$.\qquad \text{(McCoy's property).}} \]
		
		\begin{proof}
			Let $f$ and $g$ be elements of $B$ such that $fg=0$. By the amount formula in Definition \ref{amountalgebras}, there are non-negative integers $m$ and $n$ such that \[A_f ^{m+1} A_g ^{n+1}=(0). \] Since $R$ is reduced, $A_f A_g = (0)$. So, we have already proved that $B$ is Armendariz. Now let $f$ be a zero-divisor in $B$. By definition, there is a nonzero element $g$ in $B$ such that $fg=0$. Since $B$ is Armendariz $A_f A_g =(0)$. Note that $g$ is nonzero and so $A_g$ is a nonzero ideal of $R$. Take $r$ to be a nonzero element of $A_g$. Therefore, $rA_f =(0)$. This implies that $A_{rf}=(0)$. Hence, $fr=0$, i.e. McCoy's property holds. This completes the proof. 	\end{proof}
	\end{theorem}
	
	\begin{theorem}
		
		\label{primesamountextension}
		
		Let $B$ be an amount $R$-algebra. Then $P$ is a prime ideal of $R$ if and only if $P^{\epsilon}$ is a prime ideal of $B$.
		\begin{proof}
			
			Let $P$ be a prime ideal of $R$ and $fg\in P^{\epsilon}$ for arbitrary $f,g \in B$. It is clear that $A_{fg} \subseteq P$. On the other hand, by the amount formula in Definition \ref{amountalgebras}, there are non-negative integers $m$ and $n$ such that \[A_f ^m A_g ^n A_{fg} = A_f ^{m+1} A_g ^{n+1}.\] Therefore, $A_f ^{m+1} A_g ^{n+1} \subseteq P.$ Since $P$ is prime, either $A_f \subseteq P$ or $A_g \subseteq P$. This means either $f\in P^{\epsilon}$ or $g\in P^{\epsilon}$. Note that $P^{\epsilon} \neq B$. Therefore, $P^{\epsilon}$ is a prime ideal of $B$.
			
			Now let $P^{\epsilon}$ be a prime ideal of $B$ and $r$ and $s$ be elements of $R$ such that $rs\in P$. This implies that $A_{rs} = (rs) \subseteq P$. So, $rs\in P^{\epsilon}$. From this, we obtain that either $r\in P^{\epsilon}$ or $s\in P^{\epsilon}$ which is equivalent to say that either $r\in P$ or $s\in P$ and this completes the proof.
		\end{proof}	
	\end{theorem}
	
	In the following, we recall the definition of $n$-absorbing and strongly $n$-absorbing ideals, and also the definition of $\omega_R(I)$ \cite{AndersonBadawi2011}. For more on $n$-absorbing ideals and related topics refer to the recent survey paper \cite{Badawi2017}.
	
	\begin{definition} \label{n-absorbing}
		Let $R$ be a ring.\begin{enumerate}
			\item A proper ideal $I$ of $R$ is an $n$-absorbing ideal of $R$, if whenever $r_1 \cdots r_{n+1} \in I$ for $r_1, \ldots, r_{n+1} \in R$, then there are $n$ of the $r_i$'s whose product is in $I$. 
			
			\item If there is a positive integer $n$ such that $I$ is an $n$-absorbing ideal of $R$, then \[\omega_R(I)= \min \{n \colon \text{ $I$ is an $n$-absorbing ideal of $R$}\}.\] Otherwise, $\omega_R(I)=\infty.$
			
			\item A proper ideal $I$ of $R$ is a strongly $n$-absorbing ideal of $R$ if whenever $I_1 \cdots I_{n+1} \subseteq I$ for some ideals $I_1, \ldots, I_{n+1}$ of $R$, then there are $n$ of the $I_i$'s whose product is a subset of $I$.
		\end{enumerate}
	\end{definition}

	The proof of the following statement is straightforward but we bring it only for the sake of reference.
	
	\begin{proposition}
		\label{omegaextension}
		
		If $I$ is an ideal of a ring $R$, then $\omega_R (I) \le \omega_{R[X]} (I[X]) \le \omega_{R[[X]]} (I[[X]])$.
	\end{proposition}

	\begin{definition}
		\label{Gaussianamountalgebra}
		We say an amount $R$-algebra $B$ is Gaussian if $A_{fg} = A_f A_g$ for all $f,g \in B$, where $A$ is the amount function defined in Definition \ref{amonutfunctionsdef}. 
	\end{definition}
	
	\begin{proposition}
		If an amount $R$-algebra $B$ is Gaussian then it is Armendariz.
		
		\begin{proof}
			Straightforward. 
		\end{proof}
		
	\end{proposition}
	
	\begin{examples} 
		
		\label{Gaussianamount}
		
		\begin{enumerate}
			
			\item(A general example) Let $B$ be an amount $R$-algebra such that $A_f$ is a cancellation ideal of $R$ for all nonzero elements $f$ in $B$. Then $B$ is Gaussian.
			
			\item Let us recall that a ring $R$ is Gaussian if $c(fg) = c(f)c(g)$ for all $f,g\in R[X]$ \cite{Tsang1965}. Now it is clear that if $R$ is a Gaussian ring, then the amount $R$-algebra $R[X]$ is Gaussian.
			
			\item If $D$ is a Dedekind domain, then the amount $D$-algebra $D[[X]]$ is Gaussian (Use Theorem 2.6 in \cite{EpsteinShapiro2016} and this fact that each nonzero ideal of a Dedekind domain is a cancellation ideal).
		\end{enumerate}
		
	\end{examples}
	
	\begin{lemma}
		
		\label{absorbingamountextension}
		
		Let $R$ be a ring and $I$ a proper ideal of $R$. Also, let $B$ be an amount $R$-algebra. If $I^{\epsilon}$ is $n$-absorbing, then so is $I$. Moreover, $\omega_R(I) \le \omega_B(I^{\epsilon})$.
		
		\begin{proof} Let $r_1 \cdots r_{n+1} \in I$. So, $A_{r_1 \cdots r_{n+1}} = (r_1 \cdots r_{n+1}) \subseteq I$. This implies that $r_1 \cdots r_{n+1}\in I^{\epsilon}.$ Since $I^{\epsilon}$ is $n$-absorbing, $r_1 \cdots r_{i-1} r_{i+1} r_n$ is in $I^{\epsilon}$ for some index $i$. So, \[r_1 \cdots r_{i-1} r_{i+1} r_n \in I^{\epsilon} \cap R = I.\] Now, it is clear that $\omega_R(I) \le \omega_B(I^{\epsilon})$.\end{proof}
		
	\end{lemma}
	
	\begin{theorem}
		
		\label{ABconjecture1}
		
		Let $R$ be a ring such that any $n$-absorbing ideal $I$ of $R$ is strongly $n$-absorbing for any positive integer $n$. Let $B$ be an amount $R$-algebra. If $B$ is Gaussian then $\omega_B(I^{\epsilon}) = \omega_R(I)$.
		
		\begin{proof}
			By Lemma \ref{absorbingamountextension}, $\omega_R(I) \le \omega_B(I^{\epsilon})$. Let $I$ be a proper ideal of $R$ such that $\omega_R(I)=n$ for a positive integer $n$. Our claim is that $I^{\epsilon}$ is an $n$-absorbing ideal of $B$. Assume that \[f_1 \cdots f_{n+1} \in I^{\epsilon},\] for arbitrary $f_1, \ldots, f_{n+1} \in B$. 
			
			It is clear that $A_{f_1 \cdots f_{n+1}} \subseteq I$. Since $B$ Gaussian, $A_{f_1 \cdots f_{n+1}}=A_{f_1} \cdots A_{f_{n+1}}$. 
			By assumption, $I$ is a strongly $n$-absorbing ideal of $R$. 
			
			Therefore, $A_{f_1} \cdots A_{f_{i-1}}A_{f_{i+1}} \cdots A_{f_{n+1}} \subseteq I$ for some $i$. This implies that \[A_{f_1 \cdots f_{i-1} f_{i+1} \cdots f_{n+1}} \subseteq I.\] And this means that \[f_1 \cdots f_{i-1} f_{i+1} \cdots f_{n+1} \in I^{\epsilon}.\] So, we have already proved that $n=\omega_R(I) \le \omega_B(I^{\epsilon}) \le n$. Finally, it is easy to see that $\omega_B(I^{\epsilon}) = \infty$ if and only if $\omega_R(I) = \infty$, and the proof is complete.
		\end{proof}
		
	\end{theorem}
	
	\begin{corollary}
		
		\label{AB1cor}
		Let $D$ be a Pr\"{u}fer domain. If an amount $D$-algebra $B$ is Gaussian, then $\omega_B(I^{\epsilon}) = \omega_D(I)$ for each ideal $I$ of $D$. In particular, if $I$ is an ideal of a Dedekind domain $D$, then \[\omega_{D[[X]]} (I[[X]]) = \omega_{D[X]} (I[X]) =  \omega_D (I).\]
		
		\begin{proof} Since $D$ is a Pr\"{u}fer domain, any $n$-absorbing ideal of $D$ is strongly $n$-absorbing for each positive integer $n$ \cite[Corollary 6.9]{AndersonBadawi2011}. Now by Theorem \ref{ABconjecture1}, $\omega_B(I^{\epsilon}) = \omega_R(I)$. In particular, if $D$ is a Dedekind domain, by Examples \ref{Gaussianamount}, \[\omega_{D[[X]]} (I[[X]]) = \omega_{D[X]} (I[X]) =  \omega_D (I),\] and this completes the proof.\end{proof}
	\end{corollary}
	
	\begin{theorem}
		
		\label{ABconjecture2}
		
		Let $R$ be a ring such that any $n$-absorbing ideal $I$ of $R$ is strongly $n$-absorbing for any positive integer $n$. Let $B$ be an amount $R$-algebra. If $I$ is a radical ideal of $R$, then $\omega_B(I^{\epsilon}) = \omega_R(I)$.
		
		\begin{proof} 
			
			Let $f_1 \cdots f_{n+1} \in I^{\epsilon}$. Obviously, $A_{f_1 \cdots f_{n+1}} \subseteq I$. Let $g=f_2 \cdots f_{n+1}$. By the amount formula in Definition \ref{amountalgebras}, there are non-negative integers $m,n$ such that \[A_{f_1} ^m A_g ^n A_{f_1 g} = A_{f_1} ^{m+1} A_g ^{n+1}, \] and since $A_{f_1 g} \subseteq I$, we have $A_{f_1}^{m+1} A_g ^{n+1} \subseteq I$. Take $u=\max\{m,n\}$. It is easy to see that $(A_{f_1} A_g) ^{u+1}=A_{f_1}^{u+1} A_g ^{u+1} \subseteq I$. Since $I$ is a radical ideal of $R$, we have $A_{f_1} A_g \subseteq I$. 
			
			Now let $h=f_3 \cdots f_{n+1}$. It is clear that $g=f_2 h$ and by the amount formula in Definition \ref{amountalgebras}, there are non-negative integers $k,l$ such that \[A_{f_2} ^k A_h ^l A_{f_2 h} = A_{f_2} ^{k+1} A_h ^{l+1}. \] Obviously, we have the following: \[A_{f_1} A_{f_2} ^{k+1} A_h ^{l+1} = A_{f_1} A_{f_2} ^k A_h ^l A_{f_2 h} = A_{f_1} A_{f_2} ^k A_h ^l A_{g} \subseteq I.\] Similarly, since $I$ is a radical ideal of $R$, we have $A_{f_1} A_{f_2} A_h \subseteq I$. Continuing this process, we obtain that \[A_{f_1} \cdots A_{f_{n+1}} \subseteq I.\]
			
			Now if $I$ is an $n$-absorbing ideal of $R$, then according to our assumptions, $I$ is strongly $n$-absorbing. Thus, \[A_{f_1} \cdots A_{f_{i-1}} A_{f_{i+1}} \cdots A_{f_{n+1}} \subseteq I\] for some $i$.
			
			On the other hand, by Definition \ref{amonutfunctionsdef}, the amount function $A$ is submultiplicative. Therefore,  \[A_{f_1 \cdots f_{i-1} f_{i+1} \cdots f_{n+1}} \subseteq A_{f_1} \cdots A_{f_{i-1}} A_{f_{i+1}} \cdots A_{f_{n+1}}. \] This implies that $f_1 \cdots f_{i-1} f_{i+1} \cdots f_{n+1} \in I^{\epsilon}$ and so $I^{\epsilon}$ is $n$-absorbing. 
			
			Now by considering Lemma \ref{absorbingamountextension}, the rest of the proof is similar to the proof of Theorem \ref{ABconjecture1}. This completes the proof. \end{proof}
		
	\end{theorem}
	
	Let us recall that a ring $(R,+,\cdot)$ is torsion-free if $(R,+)$ is a torsion-free group \cite{BeaumontPieace1961}.
	
	\begin{corollary} \label{AB2cor1}
		Let $R$ be a torsion-free Noetherian ring and $I$ a radical ideal of $R$. Then \[\omega_{R[[X]]} (I[[X]]) = \omega_{R[X]} (I[X]) =  \omega_R (I).\]
		
		\begin{proof} Since $R$ is Noetherian, by Theorem 2.6 in \cite{EpsteinShapiro2016}, $R[[X]]$ is an amount $R$-algebra. On the other hand, since $R$ is torsion-free, by Theorem 4.2 in \cite{DaraniPuczylowski2013}, each $n$-absorbing ideal of $R$ is strongly $n$-absorbing for any positive integer $n$. By using Theorem \ref{ABconjecture2}, the proof of this corollary is complete. \end{proof}
	\end{corollary}
	
	\begin{corollary}
		
		\label{AB2cor2}
		Let $I$ be a radical ideal of a domain $D$. If either $D$ is a Pr\"{u}fer domain or $D$ is a torsion-free valuation ring, then \[\omega_{D[[X]]} (I[[X]]) = \omega_{D[X]} (I[X]) =  \omega_D (I).\]
		
		\begin{proof} If either $D$ is a Pr\"{u}fer domain or $D$ is a torsion-free valuation ring, then by the Theorem 2.8 and the proof of Corollary 2.9 in \cite{KangParkToan2018}, in each case, $D[[X]]$ is an amount $D$-algebra. Also, in each of the mentioned cases, any $n$-absorbing ideal of $D$ is strongly $n$-absorbing (see Corollary 6.9 in \cite{AndersonBadawi2011} and Theorem 4.2 in \cite{DaraniPuczylowski2013}). In view of Theorem \ref{ABconjecture2}, the proof of this corollary is complete. \end{proof}
	\end{corollary}

	\begin{conjecture}
		Let $X$ be an indeterminate over a ring $R$. For any ideal $I$ of $R$, \[\omega_{R[[X]]} (I[[X]])  =  \omega_R (I).\]
	\end{conjecture}
	
	\subsection*{Acknowledgments} This work is supported by the Golpayegan University of Technology. Our special thanks go to the Department of Engineering Science at the Golpayegan University of Technology for providing all the necessary facilities available to us for successfully conducting this research.
	
	\bibliographystyle{amsplain}

\end{document}